\documentclass[12pt]{article}

\usepackage{amssymb} %
\usepackage{amsmath} %
\usepackage{tabularx}
\usepackage{subcaption}
\usepackage{caption}
\usepackage{calc}
\usepackage{tikz} %
\usepackage{graphicx}
\usepackage{pstricks}
\usepackage{float} %
\usepackage{multicol}
\usepackage{amsthm}
\usepackage{comment} %
\usepackage{lscape}
\usepackage{lipsum}

\usetikzlibrary{decorations.markings}
\tikzstyle{vertex}=[circle, draw, inner sep=0pt, minimum size=6pt]

\newtheorem{theorem}{Theorem}[section]

\newtheorem{lemma}[theorem]{Lemma}

\newtheorem{corollary}[theorem]{Corollary}

\newenvironment{pf}           {\noindent{\bf Proof.} }%
                                {\null\hfill$\Box$\par\medskip\medskip\medskip\medskip}
\begin{document}

\title{Graceful Labellings of Various Cyclic Snakes}
\author{Ahmad H. Alkasasbeh \footnote{ahmad84@mun.ca} \hspace{0.9in} Danny Dyer \footnote{dyer@mun.ca}  \\
Department of Mathematics and Statistics \\
Memorial University of Newfoundland \\
St. John's, Newfoundland \\
A1C 5S7 Canada}
\maketitle

\begin{abstract}

In this paper, we present a new sufficiency condition to obtain a graceful labelling for every $kC_{4n}$ snake and use this condition to label every such snake for $n=1,2,\ldots,6$. Then, we extend this result to cyclic snakes where the cycles lengths vary. Also, we obtain new results on the (near) graceful labelling of cyclic snakes based on cycles of lengths $n=6, 10, 14$, completely solving the case $n=6$.
%In this paper, we present a new sufficiency condition to obtain a graceful labelling for every $kC_{4n}$-snake and use this condition to label every such snake for $n=1,2,\ldots,6$. Then, we extend this result to cyclic snakes where the cycles vary. Also, we obtain new results on the (near) graceful labelling of cyclic snakes based on cycles of length $n=6, 10, 14$, completely solving the case for $n=6$.% variable $n_{1}C_{4m_{1}}n_{2}C_{4m_{2}}\ldots n_{i}C_{4m_{i}}$-snake
\end{abstract}
%Key words: .
%%%%%%%%%%%%%%%%%%%%%%%%%%%%%%%%%%%%%%%%%%%%%%%%%%%%%%%%%%%%%%%%%%%%%%%%%%%%%%%%%%%%%%%%%%%%%%%%%%%%%%%%%%%%%%%%%%%%%%%%%%%%%%%%%%%%%%%%%%%%%%%%%%%%%%%%%%%%%%%%%%%%%%%%%%%%%%%%%%%%%%%%%%%%%%%%%%%%%%%%%%%      C6     %%%%%%%%%%%%%%%%%%%%%%%%%%%%%%%%%%%%%%%%%%%%%%%%%%%%%%%%%%%%%%%%%%%%%%%%%%%%%%%%%%%%%%
%%%%%%%%%%%%%%%%%%%%%%%%%%%%%%%%%%%%%%%%%%%%%%%%%%%%%%%%%%%%%%%%%%%%%%%%%%%%%%%%%%%%%%%%%%%%%%%%%%%%%%%%%%%%%%%%%%%%%%%%%%%%%%%%%%%%%%%%%%%%%%%%%%%%%%%%
%\section{Introduction}
%\chapter{Graceful Labellings of Various Cyclic Snakes}\label{ch:4}
%{\let\thefootnote\relax\footnote{{Submitted to the Ars Combinatoria Journal.}}}
%%%%%%%%%%%%%%%%%%%%%%%%%%%%%%%%%%%%%%%%%%%%%%%%%%%%%%%%%%%%%%%%%%%%%%%%%%%%%%%%%INTRODUCTION%%%%%%%%%%%%%%%%%%%%%%%%%%%%%%%%%%%%%%%%%%%%%%%%%%%%%%%%%%%%%%%%%%%
\section{Introduction}
Let $G=(V,E)$ be a graph with $m$ edges. Let $f$ be a labelling defined from $V(G)$ to $\{0,1,2,\dots,m\}$ and let $g$ be the induced edge labelling defined from $E(G)$ to $\{1,2,\dots,m\}$ given by $g(uv)=|f(u)-f(v)|,$ for all $uv \in E.$ The labelling $f$ is \textit{graceful} if $f$ is an injective mapping and $g$ is a bijection. If a graph $G$ has a graceful labelling, then it is graceful.

Alternatively, let $f$ be defined from $V(G)$ to $\{0,1,2,\dots,m+1\}$ and let $g$ be the induced edge labelling defined from $E(G)$ to $A$, where $A$ is $\{1,2,\dots,m-1,m\}$ or $\{1,2,\dots,m-1,m+1\}$ given by $g(uv)=|f(u)-f(v)|,$ for all $uv \in E.$ Then $f$ is \textit{near graceful} if $f$ is an injective mapping and $g$ is a bijection. If a graph $G$ has a near graceful labelling, then it is near graceful. In this paper, every near graceful labelling we find will omit the vertex label $m$ and the edge label $m$; that is the codomain of $f$ will be $\{1,2,\dots,m-1,m+1\}$ and the codomain of $g$ will be $\{1,2,\dots,m-1,m+1\}$.
%In this paper every near graceful labelling we find will be with the co-domain of $g$ is $\{1,2,\dots,m-1,m+1\}$.

A \textit{cyclic snake} is a connected graph whose block-cutpoint graph is a path and each of the blocks is isomorphic to a fixed cycle. We define $kC_{n}$ to be a cyclic snake with $k$ blocks each of which is $C_{n}$. The \textit{string} of a $kC_n$ is a sequence of integers $\left(d_1,d_2,d_3,\ldots,d_{k-2}\right)$ where $d_i$ is the distance between the $i$th and $\left(i+1\right)$th cut vertex, counting cut vertices from one end of the snake to the other. Certainly, for fixed $n$ and $k$, the string is uniquely determined by the snake and vice versa. Note that if $k=1$ or $k=2$ the snake $kC_{n}$ has no string. A $kC_{n}$ is \textit{linear} if all entries in its string are $\left\lfloor \frac{n}{2}\right\rfloor$, and it is \textit{even} if all entries in its string are even numbers.

In \cite{rosa}, Rosa showed that all cycles $C_{n}$ with $n\equiv 0$ or $3\ ($mod$\ 4)$ are graceful. Further, he introduced a necessary condition for an Eulerian graph to be graceful, namely if $G$ is a graceful Eulerian graph with $n$ edges, then $n\equiv 0$ or $3\ ($mod$\ 4)$. A $kC_{n}$ is an Eulerian graph, and hence graceful only if $kn\equiv 0$ or $3\ ($mod$\ 4)$.% is a connected graph with $k$ blocks whose block-cut-point graph is a path and each of the $k$ blocks is isomorphic to the cycle $C_{n}$.

Moulton, in \cite{mou}, proved that a graceful labelling exists for every $kC_{3}$. Barrientos, in \cite{barrientos}, proved that the cycle $C_{n}$ has a near graceful labelling if and only if $n\equiv 1$ or $2\ ($mod$\ 4)$. This paper also showed that a graceful labelling exists for every $kC_{4}$, and for particular cases of snakes for $C_{6}$, $C_{8}$, and $C_{12}$. A complete survey of graph labellings is presented in \textit{A Dynamic Survey of Graph Labelling} \cite{Gallian}.

%In 1979, Bermond in \cite{Bermond}, proved that Dutch windmills are graceful.
We define a variable snake, $n_{1}C_{m_{1}}n_{2}C_{m_{2}}\ldots n_{i}C_{m_{i}}$, to be a combination of different $n_{j}C_{m_{j}}$, where $n_{j}C_{m_{j}}$ is connected with $n_{j+1}C_{m_{j+1}}$ by identifying a vertex in the last cycle of the $n_{j}C_{m_{j}}$ with a vertex in the first cycle of the $n_{j+1}C_{m_{j+1}}$ (other than the cut vertex). The string for a variable snake is similar to the string for a $kC_n$.% The \textit{distance} between two vertices $a$ and $b$ in a graph is the number of edges in a shortest path between $a$ and $b$.

 We will represent all the cycle labellings in this paper as $n$-tuples, with the overline elements indicating the cut vertices, when necessary.

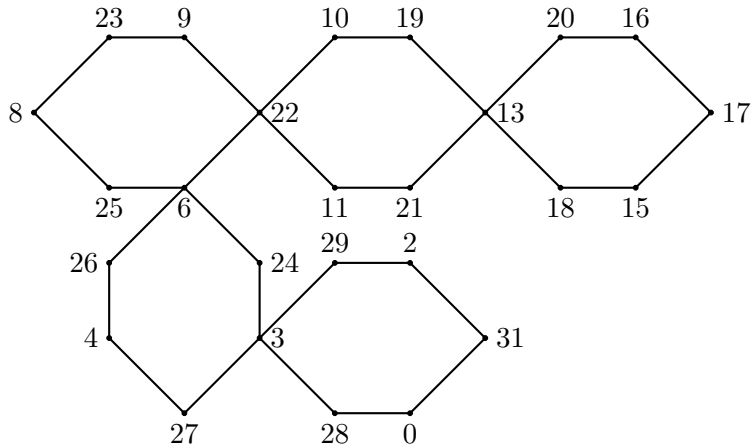
\begin{figure}[H]
\begin{center}
\scalebox{1.00}{
\begin{tikzpicture}
\draw [thick, -] (1,0) -- (2,0) -- (3,1) -- (2,2) -- (1,2) -- (0,1)  -- (1,0);
\draw[fill] (1,0) circle  [radius=.9pt];\draw[fill] (2,0) circle  [radius=.9pt];\draw[fill] (3,1) circle  [radius=.9pt];\draw[fill] (2,2) circle  [radius=.9pt];\draw[fill] (1,2) circle  [radius=.9pt];\draw[fill] (0,1) circle  [radius=.9pt];

\node [below] at (1,0) {\small $25$};
\node [below] at (2,0) {\small $6$};
\node [right] at (3,1) {\small $22$};
\node [above] at (2,2) {\small $9$};
\node [above] at (1,2) {\small $23$};
\node [left] at (0,1) {\small $8$};

%%%%%%%%%%%%%%%%%%%%%%%%%%%%%%%%%%%%% 2  %%%%%%%%%%%%%%%%%%%%%%%%%%%%%%%%%%%%%%%%%%%%%%%%%%%

\draw [thick, -] (3,1) -- (4,2) -- (5,2) -- (6,1) -- (5,0) -- (4,0)  -- (3,1);
\draw[fill] (4,0) circle  [radius=.9pt];\draw[fill] (5,0) circle  [radius=.9pt];\draw[fill] (6,1) circle  [radius=.9pt];\draw[fill] (5,2) circle  [radius=.9pt];\draw[fill] (4,2) circle  [radius=.9pt];;

\node [below] at (4,0) {\small $11$};
\node [below] at (5,0) {\small $21$};
\node [right] at (6,1) {\small $13$};
\node [above] at (5,2) {\small $19$};
\node [above] at (4,2) {\small $10$};
%\node [left] at (0,1) {\tiny $6$};

%%%%%%%%%%%%%%%%%%%%%%%%%%%%%%%%%%%  3 %%%%%%%%%%%%%%%%%%%%%%%%%%%%%%%%%%%%%%%%%%%%%%%%%%%%%

\draw [thick, -] (6,1) -- (7,2) -- (8,2) -- (9,1) -- (8,0) -- (7,0)  -- (6,1);
\draw[fill] (7,0) circle  [radius=.9pt];\draw[fill] (8,0) circle  [radius=.9pt];\draw[fill] (9,1) circle  [radius=.9pt];\draw[fill] (8,2) circle  [radius=.9pt];\draw[fill] (7,2) circle  [radius=.9pt];;

\node [below] at (7,0) {\small $18$};
\node [below] at (8,0) {\small $15$};
\node [right] at (9,1) {\small $17$};
\node [above] at (8,2) {\small $16$};
\node [above] at (7,2) {\small $20$};
%\node [left] at (0,1) {\tiny $6$};

%%%%%%%%%%%%%%%%%%%%%%%%%%%%%%%%%%%  4 %%%%%%%%%%%%%%%%%%%%%%%%%%%%%%%%%%%%%%%%%%%%%%%%%%%%%

\draw [thick, -] (2,0) -- (3,-1) -- (3,-2) -- (2,-3) -- (1,-2) -- (1,-1)  -- (2,0);
\draw[fill] (3,-1) circle  [radius=.9pt];\draw[fill] (3,-2) circle  [radius=.9pt];\draw[fill] (2,-3) circle  [radius=.9pt];\draw[fill] (1,-2) circle  [radius=.9pt];\draw[fill] (1,-1) circle  [radius=.9pt];;

\node [right] at (3,-1) {\small $24$};
\node [right] at (3,-2) {\small $3$};
\node [below] at (2,-3) {\small $27$};
\node [left] at (1,-2) {\small $4$};
\node [left] at (1,-1) {\small $26$};
%\node [left] at (0,1) {\tiny $6$};

%%%%%%%%%%%%%%%%%%%%%%%%%%%%%%%%%%%  5 %%%%%%%%%%%%%%%%%%%%%%%%%%%%%%%%%%%%%%%%%%%%%%%%%%%%%

\draw [thick, -] (3,-2) -- (4,-1) -- (5,-1) -- (6,-2) -- (5,-3) -- (4,-3)  -- (3,-2);
\draw[fill] (4,-1) circle  [radius=.9pt];\draw[fill] (5,-1) circle  [radius=.9pt];\draw[fill] (6,-2) circle  [radius=.9pt];\draw[fill] (5,-3) circle  [radius=.9pt];\draw[fill] (4,-3) circle  [radius=.9pt];;

\node [above] at (4,-1) {\small $29$};
\node [above] at (5,-1) {\small $2$};
\node [right] at (6,-2) {\small $31$};
\node [below] at (5,-3) {\small $0$};
\node [below] at (4,-3) {\small $28$};
%\node [left] at (0,1) {\tiny $6$};

\end{tikzpicture}
}
\caption{A near gracefully labelled $5C_{6}$.}\label{g1}
\end{center}
\end{figure}

In Figure~\ref{g1}, we have a near gracefully labelled $5C_{6}$ with string $\left(3,1,2\right)$. We can represent the labelling of the $5C_{6}$ in Figure~\ref{g1} in five $6$-tuples as follows: $(20$, 16, 17, 15, 18, $\overline{13})$, $(\overline{13}$, $21$, $11$, $\overline{22}$, $10$, $19)$, $(\overline{22}$, $\overline{6}$, 25, 8, 23, $9)$, $(\overline{6}$, 26, 4, 27, $\overline{3}$, $24)$, $(\overline{3}$, 29, 2, 31, 0, $28)$.
% we will discuss that later in this section.% where $n_{i}$ is the number of blocks whose cut vertex is a path and each of the $n_{i}$ blocks isomorphic to $C{m_{i}$.
% This string can be seen while starting from the cut-vertex on the top right, to left and then bottom.%we have a near gracefully labelled $5C_{6}$ snake with non linear string $3,1,2$ if you start from the cut-vertex on the top right to left then to bottom. %We call the distance between $0$ and the joint vertex by $d$, and the distance between the joint vertex and the largest value in the labelling by $\overline{d}$.
%%%%%%%%%%%%%%%%%%%%%%%%%%%%%%%%%%%%%%%%%%%%%%%%%%%%%%%%%%%%%%%%%%%%%%%%%%%%%%%%%%%%%%%%%%%%%%%%%%%%%%%%%%%%%%%%%%%%%%%%%%%%%%%%%%%%%%%%%%%%%%%%%%%%%%%%%%%%%%%%%%%%%%%
%%%%%%%%%%%%%%%%%%%%%%%%%%%%%%%%%%%%%%%%%%%%%%%%%%%%%%%  KC5  %%%%%%%%%%%%%%%%%%%%%%%%%%%%%%%%%%%%%%%%%%%%%%%%%%%%%%%%%%%%%%%%%%%%%%%%%%%%%%%%%%%%%%%%%%%%%%%%%%
In \cite{rosa}, Rosa also introduced an $\alpha$-labelling. An \textit{$\alpha$-labelling} of a graph $G$ is a graceful labelling with an extra condition which is there exists an integer $w$ such that for any edge $uv \in G$, either $f(u) \leq w < f(v)$ or $f(v) \leq w < f(u)$. Any graph with an $\alpha$-labelling is necessarily bipartite. %Hence all the graceful labelling in this paper are $\alpha$-labelling.
%Further he proved that if there exists an $\alpha$-labelling of $G$ then $G$ is a bipartite graph.

We would also like to introduce an analogue of near graceful labellings. An \textit{$\hat{\alpha}$-labelling} of a graph $G$ is a near graceful labelling with an extra condition which is there exists an integer $w$ such that for any edge $uv \in G$, either $f(u) \leq w < f(v)$ or $f(v) \leq w < f(u)$. Thus the snake in Figure~\ref{g1} has an $\hat{\alpha}$-labelling with $w=16$. In fact, by the nature of our constructions, the main results of this paper could be given in terms of $\alpha$- or $\hat{\alpha}$-labellings. While noting this to be true, we choose instead to state our results in the more familiar language of graceful and near graceful labellings.
% and we say it is $\hat{\alpha}$-labelling\\Hence, all the labellings in this paper are $\alpha$-labellings or $\hat{\alpha}$-labellings.
%%%%%%%%%%%%%%%%%%%%%%%%%%%%%%%%%%%%%%%%%%%%%%%%%%%%%%%%%%%%%%%%%%%%%%%%%%%%%%%%%%%%%%%%%%%%%%%%%%%%%%%%%%%%%%%%%%%%%%%%%%%%%%%%%%%%%%%%%%%%%%%%%%%%%%%%%%%%%%%%%%%%%%%%

In Lemmas~\ref{thirdlemma} to~\ref{secondlemma}, we give some results about (not necessarily graceful) labellings that will be useful later. In this paper, we adopt the convention that 0 is a natural number. Then, when we write $\left[x,y\right]$ with $x,y \in \mathbb{N}$ and $x<y$, we are indicating the set $\left\{z \in \mathbb{N}| x \leq z \leq y\right\}$.

If $G$ is any graph and $f$ is any labelling $G$ then we can relabel $G$ by adding a constant $c$, $h(v)=f(v)+c$. This technique preserves edge labels.
\begin{lemma}\label{thirdlemma}If $c$ is an arbitrary integer and $f$ is a labelling of a graph $G=(V,E),$ defined by $f:V(G)\rightarrow \left[0,m\right]$ then $h:V\rightarrow \left[c,c+m\right]$ defined by $h(v)=f(v)+c$, is a labelling that preserves edge labels. \end{lemma}
\begin{pf}If $v_{1}v_{2} \in E$ and $f$ is a labelling then $f(v_{1}v_{2})=\left|f(v_{1})-f(v_{2})\right|$. By definition, $h(v_{1}v_{2}) = \left|h(v_{1})-h(v_{2})\right| = \left|(f(v_{1})+c)-(f(v_{2})+c)\right|$. But then $h(v_{1}v_{2}) =  \left|f(v_{1})-f(v_{2})\right| = f(v_{1}v_{2})$. Therefore, $h$ preserves edge labels. \end{pf}
%%%%%%%%%%%%%%%%%%%%%%%%%%%%%%%%%%%%%%%%%%%%%%%%%%%%%%%%%%%%%%%%%%%%%%%%%%%%%%%%%%%%%%%%%%%%%%%%%%%%%%%%%%%%%%%%%%%%%%%%%%%%%%%%%%%%%%%%%%%%%%%%%%%%%%%%
\begin{lemma}\label{third2lemma}If $c$ is an arbitrary integer and $f$ is a labelling of a graph $G=(V,E),$ defined by $f:V(G)\rightarrow \left[0,m\right]$ then $h:V\rightarrow \left[c-m,c\right]$ defined by $h(v)=c-f(v)$, is a labelling that preserves edge labels. Further if $f$ is graceful and $c=m$, then $h(v)=m-f(v)$ is graceful.\end{lemma}%Further if $f$ is graceful and $c$ equal to the number of edges then $h$ is graceful.
\begin{pf} This proof follows the same argument as Lemma \ref{thirdlemma}.

If $f$ is graceful and $c=m$, then $h(v)=m-f(v)$ is graceful. Let $v_{1},v_{2} \in V$ such that $h(v_{1})=h(v_{2}).$ Then $m-f(v_{1})=m-f(v_{2})$ which implies $f(v_{1})=f(v_{2});$ since $f$ is injective, then $v_{1}=v_{2}.$ Therefore, $h$ is injective. Since $h$ is also edge-preserving, it is graceful.\end{pf}
%%%%%%%%%%%%%%%%%%%%%%%%%%%%%%%%%%%%%%%%%%%%%%%%%%%%%%%%%%%%%%%%%%%%%%%%%%%%%%%%%%%%%%%%%%%%%%%%%%%%%%%%%%%%%%%%%%%%%%%%%%%%%%%%%%%%%%%%%%%%%%%%%%%%%%%%%%
%In the case that we begin with a graceful labelling, we can do a little better.
%\begin{lemma}\label{firstlemma}If $f$ is a graceful labelling of a graph $G=(V,E),$ $f:V\rightarrow \left[0,m\right]$ then $h:V\rightarrow \left[0,m\right]$ defined by $h(v)=m-f(v)$ is a graceful labelling. \end{lemma}
%\begin{pf}Let $v_{1},v_{2} \in V$ such that $h(v_{1})=h(v_{2}).$ Then $m-f(v_{1})=m-f(v_{2})$ which implies $f(v_{1})=f(v_{2});$ since $f$ is injective, then $v_{1}=v_{2}.$ Therefore, $h$ is injective.
%By Lemma~\ref{newlemma}, $h$ preserves edge labels. Thus, we can conclude that $h$ is a graceful labelling.\end{pf}
%%%%%%%%%%%%%%%%%%%%%%%%%%%%%%%%%%%%%%%%%%%%%%%%%%%%%%%%%%%%%%%%%%%%%%%%%%%%%%%%%%%%%%%%%%%%%%%%%%%%%%%%%%%%%%%%%%%%%%%%%%%%%%%%%%%%%%%%%%%%%%%%%%%%%%%%
In Lemma~\ref{secondlemma} we present a similar result for near graceful labelling; we omit $1$ from the range of $f$ and the domain of $h$, because if $f(v)=1$ then $h(v)=m$, and this would contradict the definition of a near graceful labelling.
\begin{lemma}\label{secondlemma}If $f$ is a near graceful labelling of a graph $G=(V,E),$ $f:V(G)\rightarrow \{0,2,3,\dots,m-1,m+1\}$ then $h:V(G)\rightarrow \{0,2,3,\dots,m-1,m+1\}$ defined by $h(v)=(m+1)-f(v)$ is a near graceful labelling. \end{lemma}
\begin{pf} This proof follows the same argument as Lemma \ref{third2lemma}.\end{pf}
We describe the technique in Lemmas~\ref{third2lemma} and~\ref{secondlemma} as taking the complement of a (near) graceful labelling.
%The neighborhood of a vertex $v$, $\mathpzc{N}(v)$, is the set of vertices adjacent to $v:\mathpzc{N}(v) = \left\{x \in V | vx \in E\right\}.$ Let $G=(V,E)$ and $H=(V,E)$ be two graphs with $x \in V(G)$ and $y \in V(H)$. $W=(V,E)$ is a new graph obtained from $G$ and $H$ such that $V(W)=V(G) \cup V(H) \cup \left\{z\right\} \setminus \left\{x,y\right\}$ and $E(W)=E(G) \cup E(H) \cup \left\{zv: v \in \mathpzc{N}(x) \cup \mathpzc{N}(y)\right\} \setminus \left(\left\{xv: v \in \mathpzc{N}(x)\right\} \cup \left\{yu: u \in \mathpzc{N}(y)\right\}\right) $ where $z \in V(W).$

In \cite{rosa} Rosa introduced a graceful labelling for $C_{4n}$ with $n \geq 1$.
\begin{lemma}\label{rosalemma}\cite{rosa} Let $C_{4n}$ be a cycle with $4n$ edges and vertices $v_i$, for $1 \leq i \leq 4n$. Then the following labelling $f$ shows that $C_{4n}$ is graceful:
\begin{equation}
f(v_{i})=
\left\{\begin{array}{ccc}
\frac{i-1}{2} & $if$\ i\ $ is odd$, &  \\
4n+1-\frac{i}{2} & $if$\ i\ $ is even$, & i\ \leq \frac{t}{2}, \\
4n-\frac{i}{2} & $if$\ i\ $ is even$, & i > \frac{t}{2}.
\end{array} \right. \label{eq1}
\end{equation}
\end{lemma}
Barrientos in \cite{barrientos} obtained the following results.
%Barrientos proved that $kC_{4}$-snakes when $k$ is odd and even are graceful, $kC_{8}$-snakes and $kC_{12}$-snakes are graceful when $k$ is even. Also, he proved that the linear $kC_{6}$-snakes are graceful when $k$ is even and are near graceful when $k$ are odd.
\begin{theorem}\label{t1}\cite{barrientos} The $kC_{4}$ has a graceful labelling for any string. \end{theorem}
\begin{theorem}\label{t2}\cite{barrientos} The linear $kC_{6}$ is near graceful if $k$ is odd and graceful if $k$ is even. \end{theorem}
\begin{theorem}\label{t3}\cite{barrientos} The even $kC_{8}$ and $kC_{12}$ are graceful graphs. \end{theorem}
\begin{theorem}\label{t4}\cite{barrientos} The even $kC_{4n}$ with string $\left(d_{1},d_{2},\ldots,d_{k-2}\right)$, where $d_{i} \in \left\{2,4\right\},$ has a graceful labelling. \end{theorem}
\begin{theorem}\label{t5}\cite{barrientos} The even $kC_{4n}$, $4 \leq n \leq 5$, with string $\left(d_{1},d_{2},\ldots,d_{k-2}\right)$, where $d_{i} \in \left\{2,4,2n\right\}$ has a graceful labelling. \end{theorem}

In this paper, we introduce a new sufficiency condition to get a graceful labelling for every $kC_{4n}$. Then, we extend this result to $n_{1}C_{m_{1}}n_{2}C_{m_{2}}\ldots n_{i}C_{m_{i}}$. Further, we extend the results in Theorem~\ref{t2} to~\ref{t5} on (near) gracefully labelled $kC_n$ where $n = 6,8,12,16,20,24$ for all possible strings. Also, we present new results on the (near) graceful labelling of $kC_n$ where $n = 10,14$ and $k > 1$.
%%%%%%%%%%%%%%%%%%%%%%%%%%%%%%%%%%%%%%%%%%%%%%%%%%%%% Main results %%%%%%%%%%%%%%%%%%%%%%%%%%%%%%%%%%%%%%%%%%%%%%%%%%%%%%%%%%%%%%%%%%%%%%%%%%%%%%%%%%%%%%%%%%%%%%%%%%
%\section{Main Results}
\section{Graceful Labelling of $kC_{m}$ for $m\equiv  0\ ($mod$\ 4)$}\label{ss2}
%In \cite{rosa} Rosa proved that all cycles $C_{n}$ with $n\equiv 0$ or $3\ ($mod$\ 4)$ are graceful. Further, Barrientos in \cite{barrientos} proved that even $kC_{4n}$-snakes, $4 \leq n \leq 5$ with string $\left(d_{1},d_{2},\ldots,d_{k-2}\right)$, where $d_{i} \in \left\{2,4,2n\right\},$ have graceful labellings.
Since the size of a $kC_{4n}$ is $4kn\equiv  0\ ($mod$\ 4)$, we have the potential to find a graceful labelling for any $kC_{4n}$. In Theorem~\ref{tc4n}, we give a new sufficient condition which, when satisfied, shows there is a graceful labelling of a $kC_{4n}$ for any string.

%In Table~\ref{c4n4n} we see a complete useful cycle set labelling of $C_{4n}$ for $n=2,3,4,5$. The labellings
For fixed even $t$ and an arbitrary fixed positive integer $s$, an \textit{$s,t$-useful cycle with even distance $d$} is a $t$-cycle with vertices labelled from $\left[0,\frac{t}{2}-1\right] \cup \left[st-\frac{t}{2},st\right]$ and edge labels $\left[st-t+1,st\right]$, with a vertex labelled $0$ and a vertex labelled $\frac{t}{2}-1$ at distance $d$. Similarly, an \textit{$s,t$-odd useful cycle with odd distance $d$} is a $t$-cycle with vertices labelled from $\left[0, \frac{t}{2}\right] \cup \left[st-\frac{t}{2}+1,st\right]$ and edge labels $\left[st-t+1,st\right]$, with a vertex labelled $0$ and vertex labelled $st-\frac{t}{2}+1$ at distance $d$.
%A \textit{$s,t$-even useful cycle with distance $d$} is a $t$-cycle with vertices labelled from $\left[0,\frac{t}{2}-1\right] \cup \left[st-\frac{t}{2},st\right]$ and edge labels $\left[st-t+1,st\right]$, and has a vertex labelled $0$ and vertex labelled $\frac{t}{2}-1$ at distance $d$ if $d$ is even, a \textit{$s,t$-odd useful cycle with distance $d$} is a $t$-cycle with vertices labelled from $\left[0,\left\lfloor \frac{t}{2}\right\rfloor-1\right] \cup \left[st-\left\lfloor \frac{t}{2}\right\rfloor,st\right]$ and edge labels $\left[st-t+1,st\right]$, and has a vertex labelled $0$ and vertex labelled $\left(st-\left(\frac{t}{2}-1\right\right))$ at distance $d$ if $d$ is odd, where $s$ is an arbitrary fixed positive integer.

A \textit{complete $s,t$-useful cycle set} is a set of $s,t$-useful cycles of even and odd distances, the union of whose distances is $\left\{1,2,\ldots,\frac{t}{2}\right\}$. Let $C_{t}^{d}$ be an element of a complete $s,t$-useful cycle set such that the distance between the vertices labelled $0$ and $\frac{t}{2}-1$ is $d$ if $d$ is even, and the distance between the vertices labelled $0$ and $st-\frac{t}{2}+1$ is $d$ if $d$ is odd.
\begin{theorem}\label{tc4n}If there is a complete $s,4n$-useful cycle set with $s \geq 1$, then there exists a graceful labelling of any $kC_{4n}$.\end{theorem}
% With $0$ in the last cycle in any position except the cut-vertex.%can be gracefully labelled with $0$ in some position in the $(k-1)^th$ cycle with any strings ($d_{i}=1,2,\ldots,2n$ where $i=1,2,\ldots,k-3$).
\begin{pf}To prove this result, we will in fact prove a slightly more complex result: namely, that given a complete $s,4n$-useful cycle set with $s \geq 1$ and for any $k \geq 1$, then there exists a graceful labelling of any $kC_{4n}$ with $0$ in the last cycle in any position except the cut vertex.

If $0$ is in the last cycle of a $kC_{4n}$, then up to symmetry its position is uniquely determined by the distance from the last cut vertex. These distances, $d$, can only be $1,2,\ldots,2n$.

We proceed by induction on $k.$ For $k=1$ use the graceful labelling for $C_{t}$ in Lemma~\ref{rosalemma}, letting $t=4n$. The vertex labels for this graceful $1C_{4n}$ are a subset of $\left[0,4n\right]$ and the edge labels are exactly $\left[1,4n\right].$

For $k=2$, we label $2C_{4n}$ while obtaining a vertex with label $0$ at even distance $d$ from the unique cut vertex as follows. We label one cycle with the labelling used for $1C_{4n}$, with $2n-1$ added to each vertex, so that the vertex formerly labelled $0$, and now labelled $2n-1$, is the cut vertex of $2C_{4n}$. The vertices have been labelled from the set $\left[2n-1,6n-1\right]$ and, by Lemma \ref{thirdlemma}, the edge labels are $\left[1,4n\right].$ Apply the labelling $C_{4n}^{d}$ to the second cycle where $d$ is even, with the cut vertex receiving the label $2n-1$. Then this labelling of $2C_{4n}$ has all vertex labels from $\left[0,8n\right],$ and the edge labels are exactly $\left[1,8n\right],$ with no repeated vertex or edge label. That is, it is a graceful labelling of $2C_{4n}$ with the vertex labelled $0$ at even distance $d$ from the cut vertex.

In the same way we label $2C_{4n}$ while obtaining a vertex with label $0$ at odd distance $d$ from the unique cut vertex. We label one cycle with the labelling used for $1C_{4n}$, but replace each vertex label $x$ by $6n+1-x$. The vertices have been labelled from the set $\left[2n+1,6n+1\right]$ and, by Lemma \ref{thirdlemma}, the edge labels are $\left[1,4n\right].$ Apply the labelling $C_{4n}^{d}$ to the second cycle where $d$ is odd, with the cut vertex receiving the label $6n+1$. Then this labelling of $2C_{4n}$ has all vertex labels from $\left[0,8n\right],$ and the edge labels are exactly $\left[1,8n\right],$ with no repeated vertex or edge label. That is, it is a graceful labelling of $2C_{4n}$ with the vertex labelled $0$ at odd distance $d$ from the cut vertex.
%We obtain a graceful labelling of $2C_{4n}$ with $0$ at odd distance $d$ by using the previously discussed labelling, ending with the $C_{4n}^{d}$-labelling in the last cycle, then applying Lemma~\ref{firstlemma}.
	
%Let $k \geq 3$. Consider $kC_{4n}$-snake with string $d_{1},d_{2},\ldots,d_{k-2}$ with $0$ in position $d$.
Consider an arbitrary $kC_{4n}$ with $k \geq 3$, with the last entry in its string $d_{k-2}$. Let $G$ be the $(k-1)C_{4n}$ obtained by deleting a last cycle from this $kC_{4n}$. By the induction hypothesis, there is a graceful labelling of $G$ with a $0$ on the vertex distance $d_{k-2}$ from the previous cut vertex. This labelling has vertex labels that are a subset of $\left[0,4nk-4n\right]$ and the edge labels are exactly $\left[1,4kn-4n\right]$.

We label $kC_{4n}$ obtaining a vertex with label $0$ at even distance $d$ from the unique cut vertex as follows. For the first $k-1$ cycles, use the labelling of $G$ and add $2n-1$ to each vertex label, so that the final cut vertex receives label $2n-1$. Thus, the vertices have been labelled from the set $\left[2n-1,4kn-2n-1\right]$ and by Lemma \ref{thirdlemma}, the edge labels are $\left[1,4kn-4n\right]$.
 Apply the labelling $C_{4n}^{d}$ to the final cycle, with the cut vertex receiving label $2n-1$. Then this labelling of $kC_{4n}$ has all vertices labelled from $\left[0,4kn\right],$ and the edge labels are exactly $\left[1,4kn\right],$ with no repeated vertex or edge label. Thus, there is a graceful labelling of $kC_{4n}$ with the vertex labelled $0$ at even distance $d$ from the cut vertex.

We label $kC_{4n}$ obtaining a vertex with label $0$ at odd distance $d$ from the unique cut vertex as follows. For the first $k-1$ cycles, use the labelling of $G$  with a $0$ on the vertex distance $d_{k-2}$ from the previous cut vertex. Then subtract each vertex label from $4kn-2n+1$. Thus, the vertices have been labelled from the set $\left[2n+1,4kn-2n+1\right]$ and by Lemma \ref{thirdlemma}, the edge labels are $\left[1,4kn-4n\right].$
 Apply the labelling $C_{4n}^{d}$ to the final cycle, with the cut vertex receiving label $4kn-2n+1$. Then this labelling of $kC_{4n}$ has all vertices labelled from $\left[0,4kn\right],$ and the edge labels are exactly $\left[1,4kn\right],$ with no repeated vertex or edge label. Thus, there is a graceful labelling of $kC_{4n}$ with the vertex labelled $0$ at odd distance $d$ from the cut vertex.\end{pf}
%By using the previously discussed labelling, ending with the $C_{4n}^{d}$-labelling in the last cycle, and applying Lemma~\ref{firstlemma}, we obtain a graceful labelling of $kC_{4n}$ with $0$ at odd distance $d$.

In Table~\ref{c4n4n}, we give labellings for $C_{4n}^{2j}$ where $1 \leq j \leq n$ and $n \leq 6$. For each $C_{4n}^{2j}$, we can use Lemma~\ref{third2lemma} with $c=t$ to obtain $C_{4n}^{2j-1}$. Then, $\left\{C_{4n}^{m} | 1 \leq m \leq 2n\right\}$ is a complete $s,4n$-useful cycle set. Combining these sets with Theorem~\ref{tc4n}, we obtain the following corollary.
%as well as the useful cycle set taken with the complement $C_{4n}^{2j-1}$ by Lemma~\ref{newlemma} then we obtain a complete useful cycles set. Combining Theorem~\ref{tc4n} with the results of Table~\ref{c4n4n} with Lemma~\ref{newlemma} we obtain the following corollary.
\begin{corollary}If $1 \leq n \leq 6$ and $k \geq 1$ then every $kC_{4n}$ is graceful.\end{corollary}
%A \textit{complete $k,t$-useful cycle set} is a set of $k,t$-useful cycles the union of whose differences is $\left\{1,2,\ldots,\left\lfloor \frac{kt}{2}\right\rfloor\right\}$. Let $C_{4n}^{d}$ be an element of complete $k,t$-useful cycle set such that the distance between the vertices labelled $0$ and $\left\lfloor \frac{kt}{2}\right\rfloor -1$ is $d$.
%We obtain a labelling with $0$ position $j-1$ or $j+1$ by using the previously discussed labelling, ending with the $C_{4n}-$labelling in the last cycle, then apply Lemma~\ref{firstlemma}.
% The labelling of $C_{4n}^{j}$ from $\left\{1,2,3,\ldots,2n-1,4kn-2n,\ldots,4kn\right\}$ and the dege labelling is $\left\{4kn-4n+1,,\ldots,4kn-1,4kn\right\},$
%If we apply Lemma~\ref{firstlemma} on the
%\textbf{Note: If we take the complement of Rosa construction $R^{c}$ we can get $C_{4n}^{4}$ for all $n$ and its complement $R$ give us $C_{4n}^{4}$. (I have to ask Dr. Danny about it).}
\begin{table}[h!]
%\large
\begin{center}\setlength\extrarowheight{8pt}\renewcommand{\arraystretch}{1.0}
\scalebox{0.62}{
\begin{tabular}{|>{\centering\arraybackslash}p{0.1\linewidth}|>{\centering\arraybackslash}p{1.3\linewidth}|>{\centering\arraybackslash}p{0.1\linewidth}|}
\hline
            & Labelling &  \\ \hline
$C_{4}^{2}$ &$\left(0,t,\overline{1},t-2\right)$           & \cite{barrientos,Gnanajothi} \\ \hline
$C_{8}^{2}$ &$\left(0,t,\overline{3},t-4,2,t-3,1,t-1\right)$           &  \\ \hline
$C_{8}^{4}$ &$\left(0,t,1,t-1,\overline{3},t-4,2,t-3\right)$           & \cite{barrientos}  \\ \hline

$C_{12}^{2}$ &$\left(0,t,\overline{5},t-6,4,t-5,3,t-4,2,t-2,1,t-1\right)$           & \cite{barrientos} \\ \hline
$C_{12}^{4}$ &$\left(0,t,4,t-5,\overline{5},t-6,2,t-4,3,t-2,1,t-1\right)$           &  \\ \hline
$C_{12}^{6}$ &$\left(0,t,1,t-1,2,t-6,\overline{5},t-5,4,t-2,3,t-4\right)$           & \cite{barrientos} \\ \hline

$C_{16}^{2}$ &$\left(0,t,\overline{7},t-8,6,t-7,5,t-6,4,t-5,3,t-3,2,t-2,1,t-1\right)$           &  \\ \hline
$C_{16}^{4}$ &$\left(0,t,1,t-1,\overline{7},t-8,6,t-7,5,t-6,4,t-5,2,t-2,3,t-3\right)$           &  \\ \hline
$C_{16}^{6}$ &$\left(0,t,1,t-1,2,t-2,\overline{7},t-8,6,t-7,5,t-6,4,t-3,3,t-5\right)$          &  \\ \hline
$C_{16}^{8}$ &$\left(0,t,1,t-1,2,t-2,4,t-3,\overline{7},t-8,6,t-7,5,t-6,3,t-5\right)$            & \cite{barrientos} \\ \hline

$C_{20}^{2}$ &$\left(0,t,\overline{9},t-10,8,t-9,7,t-8,6,t-7,5,t-6,4,t-4,3,t-3,2,t-2,1,t-1\right)$           &  \\ \hline
$C_{20}^{4}$ &$\left(0,t,1,t-1,\overline{9},t-10,8,t-9,7,t-8,6,t-7,5,t-6,3,t-4,4,t-2,2,t-3\right)$           &  \\ \hline
$C_{20}^{6}$ &$\left(0,t,1,t-1,2,t-10,\overline{9},t-9,8,t-8,7,t-7,6,t-3,5,t-6,4,t-2,3,t-4\right)$           &  \\ \hline
$C_{20}^{8}$ &$\left(0,t,1,t-4,3,t-6,8,t-9,\overline{9},t-10,6,t-7,4,t-8,7,t-3,5,t-1,2,t-2\right)$           &  \\ \hline
$C_{20}^{10}$ &$\left(0,t,1,t-1,2,t-2,3,t-3,5,t-4,\overline{9},t-10,8,t-9,7,t-8,6,t-6,4,t-7\right)$           & \cite{barrientos}  \\ \hline

$C_{24}^{2}$ &$\left(0,t,\overline{11},t-12,10,t-11,9,t-10,8,t-9,7,t-8,6,t-7,5,t-5,4,t-4,3,t-3,2,t-2,1,t-1\right)$           &  \\ \hline
$C_{24}^{4}$ &$\left(0,t,1,t-4,\overline{11},t-12,10,t-11,9,t-10,8,t-9,7,t-7,6,t-3,5,t-5,2,t-1,3,t-8,4,t-2\right)$           &  \\ \hline
$C_{24}^{6}$ &$\left(0,t,1,t-1,2,t-12,\overline{11},t-11,10,t-10,9,t-9,8,t-8,7,t-5,6,t-7,3,t-2,4,t-3,5,t-4\right)$           &  \\ \hline
$C_{24}^{8}$ &$\left(0,t,1,t-4,2,t-1,3,t-8,\overline{11},t-12,10,t-11,9,t-7,8,t-10,7,t-5,5,t-9,4,t-3,6,t-2\right)$           &  \\ \hline
$C_{24}^{10}$ &$\left(0,t,1,t-4,2,t-1,3,t-8,4,t-5,\overline{11},t-12,10,t-11,9,t-10,8,t-9,6,t-7,7,t-3,5,t-2\right)$           &  \\ \hline
$C_{24}^{12}$ &$\left(0,t,1,t-1,2,t-2,3,t-3,4,t-5,6,t-4,\overline{11},t-12,10,t-11,9,t-10,8,t-9,7,t-7,5,t-8\right)$           &  \\ \hline
\end{tabular}}
\caption{Useful labellings $C_{4n}^{2j}$, where $t=4kn$.} \label{c4n4n}
\end{center}
\end{table}
%%%%%%%%%%%%%%%%%%%%%%%%%%%%%%%%%%%%%%%%%%%%%%%%%%%%%%%%%%%%%%%%%%%%%%%%%%%%%%%%%%%%%%%%%%%%%%%%%%%%%%%%%%%%%%%%%%%%%%%%%%%%%%%%%%%%%%%%%%%%%%%%%%%%%%%%%%%%%%%%%%%%%%%
%In \cite{rosa} Rosa introduced a graceful labelling for $C_{4n}$ with $n \geq 1$ as follows. Let $C_{4n}$ be a cycle with $4n$ edges and vertices $v_i$, for $1 \leq i \leq 4n$. Then the following labelling $f$ shows that $C_{4n}$ is graceful:
%\begin{equation}
%f(v_{i})=
%\left\{\begin{array}{ccc}
%\frac{i-1}{2} & $if$\ i\ $ is odd$, &  \\
%4n+1-\frac{i}{2} & $if$\ i\ $ is even$, & i\ \leq \frac{t}{2}, \\
%4n-\frac{i}{2} & $if$\ i\ $ is even$, & i > \frac{t}{2}.
%\end{array} \right \label{eq1}
%\end{equation}

By adding $t-4n$ to all even vertices in Equation~(\ref{eq1}) in Lemma~\ref{rosalemma} we obtain another labelling.
\begin{equation}
g(v_{i})=\left\{\begin{array}{ccc}
\frac{i-1}{2} & $if$\ i\ $is odd$, &  \\
t+1-\frac{i}{2} & $if$\ i\ $is even$, & i \leq 2n, \\ \hfill
t-\frac{i}{2} & $if$\ i\ $is even$, & i > 2n.
\end{array} \right. \label{eq2}
\end{equation}

 From Equation~(\ref{eq2}) we obtain $C_{4n}^{2}$ and $n \geq 1$, since the distance between the labels $0$ and $2n-1$ in $C_{4n}$ is $2$. The same labelling also works as $C_{4n}^{3}$, since the distance between the labels $0$ and $t-\left(2n-1\right)$ is $3$. Further, applying Lemma~\ref{third2lemma} with $c=t$ to the labelling $C_{4n}^{2}$ gives a new labelling with distance $4$ between the labels $0$ and $2n-1$, $C_{4n}^{4}$. Then, $S=\left\{C_{4n}^{m} | 2 \leq m \leq 4\right\}$ is a $s,4n$-useful cycle set (though not complete). Combining the set $S$ with Theorem~\ref{tc4n}, we obtain the following theorem.% These labellings are the only labellings we have now for $C_{4n}^{d}$ for all $n$.
\begin{theorem}\label{rrrr}The snake $kC_{4n}$ with string $\left(d_1,d_2,\ldots,d_{k-2}\right)$ has a graceful labelling if $d_{i} \in \left\{2,3,4\right\}$ for all $i$.\end{theorem}
%Let $C_{t}^{2}$ be a complete $s,t$-useful cycle set such that the distance between the vertices labelled $0$ and joint vertex ($2n-1$) is $d$ which is a useful  and the complement of that labelling gives the distance $4$.From Rosa labelling we have useful cycles for all $n$ with $d=2,4$ ($S=\left\{C_{4n}^{2,4}\right\}$ is useful cycles)and Theorem~\ref{tc4n} we can conclude Theorem~\ref{rrrr}.
%\begin{theorem}\label{rrrr} If there is a complete $s,4n$-useful cycle set with $s \geq 1$, then there a graceful labelling of any $kC_{4n}$-snake.\end{theorem
%Thus . in our labelling method in the proof of Theorem~\ref{tc4n}. (\textcolor[rgb]{1,0,0}{for d=3 we have to take the complement on the snake as I understand but here I do not know how to describe it here so for the snake we can get 2 or 3 but for cycle 2 or 4. }).This labelling is the only labelling we have now for $C_{4n}$
%%%%%%%%%%%%%%%%%%%%%%%%%%%%%%%%%%%%%%%%%%%%%%%%%%%%%%%%%%%%%%%%%%%%%%%%%%%%%%%%%%%%%%%%%%%%%%%%%%%%%%%%%%%%%%%%%%%%%%%%%%%%%%%%%%%%%%%%%%%%%%%%%%%%%%%%%%%%%%%%%%%%%%%
From the previous discussion we can gracefully label a variable snake made from any $kC_{4n}$ with $n \leq 6$. As an example of a variable snake, consider $3C_8$, gracefully labelled via Theorem~\ref{tc4n}, with the vertices labelled from the set $\left[0,24\right]$ and the edge labels $\left[1,24\right]$. Then form a new labelling via Lemma~\ref{thirdlemma} with $c=5$, so that a vertex in the last cycle obtains the label $5$. Then add any $C_{12}^{2j}$ (from Table~\ref{c4n4n}) to $3C_8$. We obtain a gracefully labelled $3C_{8}1C_{12}$. More generally, if we have a complete $k,4i$-useful cycle set for all $1 \leq i \leq n$, then, by combining these sets with Theorem~\ref{tc4n}, we obtain the following corollary.
\begin{corollary}\label{tc4nn}If there is a complete $s,4i$-useful cycle set with $s \geq 1$ for all $1 \leq i \leq n$, $j \geq 1$, and $1 \leq m_{1},m_{2},\ldots,m_{j} \leq n$, with $n_{1},n_{2},\ldots,n_{j}$ positive integers then every $n_{1}C_{4m_{1}}n_{2}C_{4m_{2}}\ldots n_{j}C_{4m_{j}}$ is graceful.\end{corollary}
%\textbf{\textcolor[rgb]{1,0,0}{Dr. Danny: What do you think about the following corollary ?\ I need to discuss the statment and the proof in our monday meeting}}
%\begin{corollary}\label{tc4nnn}If there exists a graceful labelling of any $n_{t}C_{4m_{t}}$-snake with $1 \leq t \leq j$ then every $n_{1}C_{4m_{1}}n_{2}C_{4m_{2}}\ldots n_{j}C_{4m_{j}}$-snake is graceful.\end{corollary}
\begin{pf}We follow the same method as in the proof of Theorem~\ref{tc4n}. That is, we will proceed by induction on the number of cycles $k=n_{1}+n_{2}+\cdots+n_{j}$, and will prove the slightly more complex result that if there is a complete $s,4i$-useful cycle set with $s \geq 1$ for all $1 \leq i \leq n$, $j \geq 1$, and $1 \leq m_{1},m_{2},\ldots,m_{j} \leq n$, with positive integers $n_{1},n_{2},\ldots,n_{j}$ then every $n_{1}C_{4m_{1}}n_{2}C_{4m_{2}}\ldots n_{j}C_{4m_{j}}$ is graceful with the label $0$ in any position in the last cycle except the cut vertex. For $k=1$ we have a $k,4i$-useful cycle set from Theorem~\ref{tc4n}, so there exists a graceful labelling of any $1C_{4n}$ with $0$ in any position in the last cycle.

Consider an arbitrary $n_{1}C_{4m_{1}}n_{2}C_{4m_{2}}\ldots n_{j}C_{4m_{j}}$ with $k > 1$ and $m=\sum\limits_{l=1}^{j} n_{l}\left(4m_{l}\right)$, the total number of edges. Let $G$ be the graph obtained by deleting a last $4m_{j}$ cycle from $n_{1}C_{4m_{1}}n_{2}C_{4m_{2}}\ldots n_{j}C_{4m_{j}}$. By the induction hypothesis, there is a graceful labelling of $G$ with a $0$ on the vertex distance $d_{k-2}$ from the previous cut vertex. This labelling has vertex labels that are a subset of $\left[0,m-4m_{j}\right]$ and the edge labels are exactly $\left[1,m-4m_{j}\right]$.

We label $n_{1}C_{4m_{1}}n_{2}C_{4m_{2}}\ldots n_{j}C_{4m_{j}}$ obtaining a vertex with label $0$ at even distance $d$ from the last cut vertex as follows. For the first $k-1$ cycles, use the labelling of $G$ and add $2m_{j}-1$ to each vertex label, so that the final cut vertex receives label $2m_{j}-1$. Thus, the vertices have been labelled from the set $\left[2m_{j}-1,m+2m_{j}-1\right]$ and by Lemma \ref{thirdlemma}, the edge labels are $\left[1,m-4m_{j}\right]$.

Now the labelling of $G$ has the label $2m_{j}-1$ at the last cycle. So apply the labelling of $C_{4m_{j}}^{d}$ to the last cycle of $G$, with the cut vertex receiving label $2m_{j}-1$ and $0$ at position $d$ (even distance) from the cut vertex. Then this labelling has all vertices labelled from $\left[0,m\right],$ and edge labels exactly $\left[1,m\right]$, with no repeated vertex or edge label. Thus, there is a graceful labelling of $n_{1}C_{4m_{1}}n_{2}C_{4m_{2}}\ldots n_{i}C_{4m_{j}}$ with the vertex labelled $0$ at even distance $d$ from the cut vertex.% where $t=\sum\limits_{l=1}^{j} n_{l}\left(m_{l}\right)$

We label $n_{1}C_{4m_{1}}n_{2}C_{4m_{2}}\ldots n_{j}C_{4m_{j}}$ to obtain a vertex with label $0$ at odd distance $d$ from the last cut vertex as follows. For the first $k-1$ cycles, use the labelling of $G$ and subtract each vertex label from $m-2m_{j}+1$, so that the final cut vertex receives label $m-2m_{j}+1$. Thus, the vertices have been labelled from the set $\left[2m_{j}-1,m-2m_{j}+1\right]$ and by Lemma \ref{thirdlemma}, the edge labels are $\left[1,m-2m_{j}+1\right]$.

Now the labelling of $G$ has the label $m-2m_{j}+1$ at the last cycle. So apply the labelling of $C_{4m_{j}}^{d}$ to the last cycle of $G$, with the cut vertex receiving label $m-2m_{j}+1$ and $0$ at position $d$ (odd distance) from the cut vertex. Then this labelling has all vertices labelled from $\left[0,m\right],$ and edge labels exactly $\left[1,m\right]$, with no repeated vertex or edge label. Thus, there is a graceful labelling of $n_{1}C_{4m_{1}}n_{2}C_{4m_{2}}\ldots,n_{i}C_{4m_{j}}$ with the vertex labelled $0$ at odd distance $d$ from the cut vertex.% where $t=\sum\limits_{l=1}^{j} n_{l}\left(m_{l}\right)$
\end{pf}
\section{Graceful Labelling of $kC_{m}$ for $m\equiv  2\ ($mod$\ 4)$}\label{ss3}
% (\textcolor[rgb]{1,0,0}{Dr. Danny, I did my best for the following paragraph. I know you will dislike it I repeat myself but for sure I need your feedback}).
In Section~\ref{ss2}, we proved that if there is a complete $s,4n$-useful cycle set with $s \geq 1$, then there exists a graceful labelling of any $kC_{4n}$ for $k \geq 1$. The next natural question is can we prove the same results for $kC_{4n+2}$?

%but for $kC_{4n+2}$-snakes it is harder, we know for $n=1,2$ if $k$ is even its graceful and if $k$ is odd its near graceful so we can guess it is similar for any $n$. The problem here is if $k$ is even then $k+1$ is odd or if $k$ is odd then $k+1$ is even. This is why for a graceful labelling, we can do better.
In $kC_{4n}$ we have the nice property that $kC_{4n}$ is always graceful regardless of the parity of $k$. For $kC_{4n+2}$ we will obtain graceful or near graceful labellings depending on the value of $k$, because the size of a $kC_{4n+2}$ is $4kn +2k$ which is congruent to $0$ modulus $4$ for $k$ even and $4kn +2k$ which is congruent to $2$ modulus $4$ for $k$ odd. Thus, we would essentially need to find two complete useful cycle sets for $4n+2$ because we are trying to change a graceful labelling to a near graceful one, or the reverse.%because if $k$ is even then $k+1$ is odd or if $k$ is odd then $k+1$ is even, thus we need one set for graceful and another one for near graceful. %Also, it is not easy to find a complete cycle sets for $4n+2$.

For $kC_{4n}$ we obtained the complete cycle set by taking the complement as in Lemmas~\ref{third2lemma} to~\ref{secondlemma}. Here for $kC_{4n+2}$ we need to omit $1$ from the labelling of the the near graceful cycles, because if we use $1$ and take the complement we will obtain a labelling that is not near graceful. For example, if we take $\left(0,7,3,1,2,5\right)$ as a labelling of $1C_{6}$, the complement would have $6$ in the resulting vertex labelling and hence would not be near graceful. In this section we prove that there exists a (near) graceful labelling of any $kC_{6}$ (Theorem~\ref{ac6}), because we found an analogue of a complete $k,6$-useful cycle set in Table~\ref{c66}, and did not use $1$ for the near graceful useful cycles.
An exhaustive analysis shows that no labelling of $C_{10}$ exists that uses the labels $\left[0,5\right]$ and $\left[t-5,t+1\right] \setminus \left\{t\right\}$ that omits the label $1$. Thus, despite their effectiveness in the $4n$-cycle case, complete cycle sets cannot help us label even $kC_{10}$.
% Unfortunately, for $1C_{10}$ we can't find any complete useful cycle set because $1$ must show up as a label for the near graceful useful cycles we can check that based on exhaustive analysis of all the cases, analyzing all the cases is not complicated. Then we can not show that there exists a (near) graceful labelling of any $kC_{10}$-snake.
%We note that from the results of Theorem~\ref{ac6} we can build any $kC_{6}$-snake with $0$ at any position except the cut-vertex, by using Table~\ref{c66} which has a complete useful $k\left(6\right)$-cycle set. So, for $n\equiv  2\ ($mod$\ 4)$ if we have a complete useful $k\left(4n+2\right)$-cycle set with $k \geq 2$ and $k$ is even, then there exists a graceful labelling of any $kC_{4n+2}$-snake with $0$ in the last cycle in any position except the cut-vertex. Further, if we have a complete useful $k\left(4n+2\right)$-cycle set with $k \geq 1$ and $k$ is odd, then there exists a near graceful labelling of any $kC_{4n+2}$-snake with $0$ in the last cycle in any position except the cut-vertex.

In \cite{barrientos}, Barrientos proved that the $kC_{4}$ has a graceful labelling with any string, as summarized in Theorem~\ref{t1}. Recall the result of Barrientos from Theorem~\ref{t2}: that the linear $kC_{6}$ is near graceful if $k$ is odd and graceful if $k$ is even. In Theorem~\ref{ac6} we prove (near) graceful labellings exist for any $kC_{6}$.

In Table~\ref{c66} we see four labellings of $C_{6}.$ The labellings $C_{6}^{a}$ and $C_{6}^{b}$ use edge labels $\left[6k-5,6k+1\right] \setminus \left\{6k\right\}$. The labellings $C_{6}^{c}$ and $C_{6}^{d}$ use edge labels $\left[6k-6,6k\right] \setminus \left\{6k-5\right\}$.
\begin{table}[H]
\large
\begin{center}%\setlength\extrarowheight{7pt}\renewcommand{\arraystretch}{1.0}
\scalebox{1.00}{
\begin{tabular}{|c|c|c|}
\hline
 & Labelling    \\
\hline
$C_{6}^{a}$ &\ $\left(0,6k+1,2,6k-1,\overline{3},6k-2\right)$     \\ % It was A
\hline
$C_{6}^{b}$ &\ $\left(0,6k+1,\overline{3},6k-2,2,6k-1\right)$       \\ % It was C
\hline
$C_{6}^{c}$ &\ $\left(0,6k,\overline{3},6k-3,1,6k-1\right)$       \\ % It was B
\hline
$C_{6}^{d}$ &\ $\left(0,6k,1,6k-1,\overline{3},6k-3\right)$       \\ % It was D
\hline
\end{tabular}}
\caption{Some useful labellings of $C_{6}$.} \label{c66}
\end{center}
\end{table}
\begin{theorem}\label{ac6}If $k \geq 1$ then there exists a (near) graceful labelling of any $kC_{6}$.\end{theorem}% with , for any strings with $0$ in the last cycle in any position except the cut-vertex.
\begin{pf}As in proof of Theorem~\ref{tc4n}, we prove a slightly more complex result. Namely, we prove that if $k \geq 1$, then there exists a (near) graceful labelling of any $kC_{6}$ with $0$ in the last cycle in any position except the cut vertex.

If $0$ is in the last cycle of a $kC_{6},$ then up to symmetry its position is uniquely determined by the distance from the last cut vertex. These distances, $d$, can only be $1,2$ or $3.$

We proceed by induction on $k.$ For $k=1,$ use $C_{6}^{a}$ or $C_{6}^{b}$ in Table~\ref{c66} with $k=1$ which will make it near graceful.
For $k=2,$ use $\left(4,7,2,\overline{9},5,6\right)$ and $\left(\overline{9},0,12,1,11,3\right)$ to obtain a labelling with $d=1$; $\left(8,5,10,\overline{3},7,6\right)$ and $\left(\overline{3},12,0,11,1,9\right)$ to obtain a labelling with $d=2$; and $\left(4,7,2,\overline{9},5,6\right)$ and $\left(9,1,11,0,12,3\right)$ to obtain a labelling with $d=3$.

%Case 1:Consider $kC_{6}$-snake with $k \geq 4$ and $k$ is even. Consider the labelling use to form $(k-1)C_{6}$, by the induction hypothesis we know the near graceful labelling exists of the $(k-1)C_{6}$ with the vertex labels are a subset of $\left[0,6k-5\right] \setminus \left\{6k-6\right\}$ and the edge labels are exactly $\left[1,6k-5\right] \setminus \left\{6k-6\right\}$ with $0$ in the last cycle in any position except the cut-vertex.
\textbf{Case 1}: Consider an arbitrary $kC_{6}$ with $k \geq 4$ and $k$ even, with the last entry in the string $d_{k-2}$. Let $G$ be the $(k-1)C_{6}$ obtained by deleting a last cycle from this $kC_{6}$. By the induction hypothesis, there is a near graceful labelling of $G$ with a $0$ on the vertex distance $d_{k-2}$ from the previous cut vertex. This labelling has vertex labels that are a subset of $\left[0,6k-5\right] \setminus \left\{6k-6\right\}$ and the edge labels are exactly $\left[1,6k-5\right] \setminus \left\{6k-6\right\}$.% with $0$ in the last cycle in any position except the cut-vertex.
%Then in any $(k-1)C_{6}$ the string is $d_{1},d_{2},\ldots,d_{k-3}$ with $d=k-2$. Then any $(k-1)C_{6}$ can be near gracefully labelled with $0$ in some position in the $(k-1)^th$ cycle. The vertex labels for this near graceful $(k-1)C_{6}$ are a subset of $\left[0,6k-5\right] \setminus \left\{6k-6\right\}$ and the edge labels are exactly $\left[1,6k-5\right] \setminus \left\{6k-6\right\}.$
%We label $kC_{6}$ as follows. The first $(k-1)$th cycle, for a $(k-1)C_{6},$ we label them by adding $3$ to each vertex with the final cut-vertex getting label $3.$ The vertices have been labeled from the following set: $\left[3,6k-2\right] \setminus \left\{6k-3\right\}$ and by Lemma \ref{thirdlemma}, the edge labels are $\left[1,6k-5\right] \setminus \left\{6k-6\right\}.$
%Apply the labelling of $C_{6}^{b}$ or $C_{6}^{d}$ to the final cycle. Then this labelling has all verticies labelled from $\left[0,6k\right],$ and edge labels are exactly $\left[1,6k\right],$ with no repeated vertex or edge label. Then there is a graceful labelling of $kC_{6}$ and $0$ in the $d=2$ position.
%We obtain a labelling with $0$ in the $d=1$ position by using the previously discussed labelling, ending with the $C_{6}^{b}$-labelling in the last cycle, then apply Lemma~\ref{firstlemma}. We obtain $d=3$ by using the $C_{6}^{d}$-labelling in the last cycle, and apply Lemma~\ref{firstlemma}.

We label $kC_{6}$ obtaining a vertex with label $0$ at distance $d$ from the cut vertex as follows. For the first $k-1$ cycles, use the labelling of $G$ and add $3$ to each vertex, so that the final cut vertex receives label $3$. Thus, the vertices have been labelled from the set $\left[3,6k-2\right] \setminus \left\{6k-3\right\}$ and by Lemma \ref{thirdlemma}, the edge labels are $\left[1,6k-5\right] \setminus \left\{6k-6\right\}.$

 Apply the labelling $C_{6}^{c}$ or $C_{6}^{d}$ to the final cycle, with the cut vertex receiving label $3$. Then this labelling of $kC_{6}$ has all vertices labelled from $\left[0,6k\right],$ and the edge labels are exactly $\left[1,6k\right],$ with no repeated vertex or edge label. Then there is a graceful labelling of $kC_{6}$ with $0$ in the $d=2$ position.

We obtain a labelling with $0$ in the $d=1$ position by using the previously discussed labelling, ending with the $C_{6}^{c}$-labelling in the last cycle, then applying Lemma~\ref{third2lemma}. We obtain $d=3$ by using the $C_{6}^{d}$-labelling in the last cycle, and applying Lemma~\ref{third2lemma}.

\textbf{Case 2}: Consider an arbitrary $kC_{6}$ with $k \geq 3$ and $k$ odd, with the last entry in the string $d_{k-2}$. We proceed in the same fashion as in Case 1, labelling all vertices except those in the final cycle, with vertices labelled from the set $\left[3,6k-3\right]$ and edges labelled $\left[1,6k-6\right]$. %Let $G$ be the $(k-1)C_{6}$ obtained by deleting a last cycle from this $kC_{6}$. By the induction hypothesis, there is a graceful labelling of $G$ with a $0$ on the vertex distance $d_{k-2}$ from the previous cut vertex. This labelling has vertex labels that are a subset of $\left[0,6k-6\right]$ and the edge labels are exactly $\left[1,6k-6\right]$.% with $0$ in the last cycle in any position except the cut-vertex.
%We label $kC_{6}$ obtaining a vertex with label $0$ at distance $d$ from the unique cut vertex as follows. For the first $k-1$ cycles, use the labelling of $G$ and add $3$ to each vertex, so that the final cut vertex receives label $3$. Thus, the vertices have been labeled from the set$\left[3,6k-3\right]$ and by Lemma \ref{thirdlemma}, the edge labels are $\left[1,6k-6\right]$.

 Apply the labelling $C_{6}^{a}$ or $C_{6}^{b}$ to the final cycle, with the cut vertex receiving label $3$. Then this labelling of $kC_{6}$ has all vertices labelled from $\left[0,6k+1\right] \setminus \left\{1,6k\right\}$, and the edge labels are exactly $\left[1,6k+1\right] \setminus \left\{6k\right\}$, with no repeated vertex or edge label. Then there is a graceful labelling of $kC_{6}$ with $0$ in the $d=2$ position.

We obtain a labelling with $0$ in the $d=1$ position by using the previously discussed labelling, ending with the $C_{6}^{b}$-labelling in the last cycle, then applying Lemma~\ref{secondlemma}. We obtain $d=3$ by using the $C_{6}^{a}$-labelling in the last cycle, and applying Lemma~\ref{secondlemma}.
\end{pf}

In Theorem~\ref{c10new} we prove (near) graceful labellings exist for some $kC_{10}$.
In Table~\ref{c10c10} we see $7$ labellings of $C_{10}.$ The labellings of $C_{10}^{a},C_{10}^{e},C_{10}^{g}$, and $C_{10}^{h}$ use edge labels $\left[10k-9,10k+1\right] \setminus \left\{10k\right\}$. The labellings of $C_{10}^{b},C_{10}^{c},C_{10}^{d}$, and $C_{10}^{f}$ use edge labels $\left[10k-10,10k\right] \setminus \left\{10k-9\right\}$. Since we cannot rely on the uniformity of a complete cycle set, this theorem uses a variety of different techniques to achieve similar effects.
%\begin{table}[H]
%\begin{center}\setlength\extrarowheight{7pt}\renewcommand{\arraystretch}{1.0}
%\scalebox{0.75}{
%\begin{tabular}{|c|c|c|}
%\hline
\begin{table}[h!]
\large
\begin{center}\setlength\extrarowheight{8pt}\renewcommand{\arraystretch}{1.0}
\scalebox{0.65}{
\begin{tabular}{|>{\centering\arraybackslash}p{0.1\linewidth}|>{\centering\arraybackslash}p{1.3\linewidth}|>{\centering\arraybackslash}p{0.1\linewidth}|}
\hline
 & Labelling    \\
\hline
$C_{10}^{a}$ &\ $\left(0,10k+1,4,10k-2,3,10k-4,\overline{5},10k-3,1,10k-1\right)$    \\  % for d=5,4 odd     %it was E
\hline
$C_{10}^{b}$ &\ $\left(0,10k,1,10k-7,3,10k-3,\overline{4},10k-1,2,10k-2\right)$       \\ % for d=5  even   %it was F
\hline
$C_{10}^{c}$ &\ $\left(0,10k,2,10k-4,\overline{4},10k-6,1,10k-2,3,10k-1\right)$       \\ % for d=4 even
\hline
$C_{10}^{d}$ &\ $\left(0,10k+1,3,10k-2,1,10k-3,4,10k-4,\overline{5},10k-1\right)$       \\ % for d=3 odd
\hline
$C_{10}^{e}$ &\ $\left(0,10k,2,10k-2,1,10k-7,3,10k-3,\overline{4},10k-1\right)$       \\ % for d=3 even
\hline
$C_{10}^{f}$ &\ $\left(0,10k+1,\overline{5},10k-4,4,10k-3,3,10k-2,1,10k-1\right)$       \\ % for d=2,1 odd
\hline
$C_{10}^{g}$ &\ $\left(0,10k,\overline{4},10k-6,2,10k-4,3,10k-2,1,10k-1\right)$       \\ % for d=2,1 even
%\hline
%$C_{10}^{h}$ &\ $\left(0,10k+1,2,10k-1,1,10k-3,3,10k-4,4,10k-5\right)$       \\ % for d=2,1 odd
\hline
\end{tabular}}
\caption{Some useful labellings of $C_{10}$.} \label{c10c10}
\end{center}
\end{table}
%which statement is better:
\begin{comment}
\begin{table}[H]
\begin{center}\setlength\extrarowheight{7pt}\renewcommand{\arraystretch}{1.0}
\scalebox{0.65}{
\begin{tabular}{|l|l|l|}
\hline
d & $k$ is odd                       & $k$ is even                        \\ \hline
5 & $C_{10}^{a}$-subtracting $10k-5$ & $C_{10}^{b}$-subtracting $10k-5$   \\ \hline
4 & $C_{10}^{a}$-adding $5$          &$C_{10}^{c}$-adding $4$             \\ \hline
3 & $C_{10}^{e}$-subtracting $10k-5$  & $C_{10}^{d}$-subtracting $10k-5$   \\ \hline
2 & $C_{10}^{g}$-adding $5$          &$C_{10}^{f}$-adding $4$             \\ \hline
1 &                                        &                               \\ \hline
\end{tabular}}
\caption{Manual for $C_{10}$}
\end{center}
\end{table}
\end{comment}
\begin{theorem}\label{c10new}
The $kC_{10}$ ($k \geq 1$) with string ($d_{1},d_{2},\ldots,d_{k-2}$) is graceful if $k$ is even and near graceful if $k$ is odd and one of the following is true:
\begin{enumerate}
	\item  $d_{i} \in \left\{4,5\right\}$ if $i$ is odd and $d_{i}=5$ if $i$ is even, where $1 \leq i \leq k-2$,
	\item  $d_{i} \in \left\{3,4\right\}$ if $i$ is odd and $d_{i}=4$ if $i$ is even, where $1 \leq i \leq k-2$,
	\item  $d_{i} \in \left\{2,3\right\}$ if $i$ is odd and $d_{i}=3$ if $i$ is even, where $1 \leq i \leq k-2$, or
	\item  $d_{i} \in \left\{1,2\right\}$ if $i$ is odd and $d_{i}=2$ if $i$ is even, where $1 \leq i \leq k-2$.
\end{enumerate}\end{theorem}
\begin{pf}The proof is similar to the proof of Theorem~\ref{ac6} with some changes to the relabelling technique we use on vertex and edge labels.

\textbf{Case 1}: $d_{i} \in \left\{4,5\right\}$ if $i$ is odd and $d_{i}=5$ if $i$ even.

As in the proof of Theorem~\ref{ac6}, we prove that the $kC_{10}$ ($k \geq 1$) with string ($d_{1},d_{2},\ldots,d_{k-2}$) satisfying the condition of part one with $0$ in the $d=4$ or $d=5$ position in the last cycle and near graceful if $k$ is odd with $10k+1$ in the $d=5$ position in the last cycle.% If $0$ is in the last cycle of a $kC_{10},$ then up to symmetry its position is uniquely determined by the distance from the last cut-vertex. These distances, $d$, can only be $1,2,3,4$ or $5.$

We proceed by induction on $k.$ For $k=1,$ use the labelling of $C_{10}^{a}$ in Table~\ref{c10c10} with $k=1$ which will make it near graceful.
For $k=2,$ use $\left(11,8,13,9,\overline{16},5,14,6,12,10\right)$ and $\left(0,20,2,18,1,\overline{16},3,17,7,19\right)$ to get $d=5$, and take this labelling with Lemma~\ref{third2lemma} to get $d=4$.
%%%%%%%%%%%%%%%%%%%%%%%%%%%%%%%%%%%%%%%%%%%%%%%%%%%%%%%%%%%%%%%%%%%%%%%%%%%%%%%%%%%%%%%%%%%%%%%%%%%%%%%%%%%%%%%%%%%%%%%%%%%%%%%%%%%%%%%%%%%%%%%%%%%%%%%%%%%%%%%%%%%%%%%%%%%%%%%%%%%%%%%%%%%%%%%%%%%%%%%%%%%     Case 1 a    %%%%%%%%%%%%%%%%%%%%%%%%%%%%%%%%%%%%%%%%%%%%%%%%%%%%%%%%%%%%%%%%%%%%%%%%%%%%%%%%%%%%%%
%%%%%%%%%%%%%%%%%%%%%%%%%%%%%%%%%%%%%%%%%%%%%%%%%%%%%%%%%%%%%%%%%%%%%%%%%%%%%%%%%%%%%%%%%%%%%%%%%%%%%%%%%%%%%%%%%%%%%%%%%%%%%%%%%%%%%%%%%%%%%%%%%%%%%%%%

\textbf{Case 1a}: The proof is similar to the proof of Case 1 in the proof of Theorem~\ref{ac6}. Consider $kC_{10}$ to be an arbitrary snake with a string as indicated in the condition of Case 1 with $k \geq 4$ and $k$ even. The labelling of $(k-1)C_{10}$ by the induction hypothesis has vertex labels that are a subset of $\left[0,10k-9\right]\setminus \left\{2,10k-10\right\}$ and the edge labels are exactly $\left[1,10k-9\right] \setminus \left\{10k-10\right\}$ with $10k-9$ in the vertex distance $d_{k-2}$ from the previous cut vertex.

We label $kC_{10}$ as follows. For the first $k-1$ cycles, use the labelling of $\left(k-1\right)C_{10}$ obtained by induction and then subtract each vertex label from $10k-5$. Thus, the vertices have been labelled from the set $\left[4,10k-5\right] \setminus \left\{5,10k-7\right\}$ and by Lemma \ref{third2lemma}, the edge labels are $\left[1,10k-9\right] \setminus \left\{10k-10\right\}$.

 Apply the labelling $C_{10}^{b}$ to the final cycle, with the cut vertex receiving label $4$. Then this labelling of $kC_{10}$ has all vertices labelled from $\left[0,10k\right],$ and the edge labels are exactly $\left[1,10k\right],$ with no repeated vertex or edge label. By induction, a graceful labelling of $kC_{10}$ exists, with $0$ in the $d=4$ position. (Note that a possible conflict occurs as $\left(10k-5\right)-2=10k-7$, however, in the labelling of the $(k-1)C_{10}$ no vertex is labelled $2$, therefore we can use $C_{10}^{b}$ without any restriction).

We obtain a labelling with $0$ in the $d=5$ position by using the previously discussed labelling, ending with the $C_{10}^{b}$-labelling in the last cycle, then applying Lemma~\ref{third2lemma}.

\textbf{Case 1b}: Consider $kC_{10}$  to be an arbitrary snake with a string as indicated in the condition of Case 1 with $k \geq 3$ and $k$ odd. The proof is similar to the proof of Case 1a, but instead of subtracting each vertex label from $10k-5$, add $5$ to each vertex and applying the labelling $C_{10}^{a}$ to the final cycle, with the cut vertex receiving label $5$. Then this labelling of $kC_{10}$ has all vertices labelled from $\left[0,10k\right] \setminus \left\{2\right\}$ and the edge labels are exactly $\left[1,10k+1\right] \setminus \left\{10k\right\}$, with no repeated vertex or edge label. Then we obtain a near graceful labelling of $kC_{10}$ with $10k+1$ in the $d=5$ position from the cut vertex.

\textbf{Case 2}: $d_{i} \in \left\{3,4\right\}$ if $i$ is odd and $d_{i}=4$ if $i$ is even.

As in the proof of Case 1, we prove that $kC_{10}$ ($k \geq 1$) with string ($d_{1},d_{2},\ldots,d_{k-2}$) satisfying the condition of part two with $0$ in the $d=3$ or $d=4$ position in the last cycle and near graceful if $k$ is odd with $0$ in the $d=4$ position in the last cycle.

We proceed by induction on $k.$ For $k=1,$ use the labelling of $C_{10}^{a}$ in Table~\ref{c10c10} with $k=1$ which will make it near graceful. For $k=2,$ use $\left(10,9,11,5,13,\overline{4},15,8,12,7\right)$ and $\left(0,20,2,16,\overline{4},14,1,18,3,19\right)$ to get $d=4$, and take this labelling with Lemma~\ref{third2lemma} to get $d=3$.
% by Theorem~\ref{oddcycle}.
%%%%%%%%%%%%%%%%%%%%%%%%%%%%%%%%%%%%%%%%%%%%%%%%%%%%%%%%%%%%%%%%%%%%%%%%%%%%%%%%%%%%%%%%%%%%%%%%%%%%%%%%%%%%%%%%%%%%%%%%%%%%%%%%%%%%%%%%%%%%%%%%%%%%%%%%%%%%%%%%%%%%%%%%%%%%%%%%%%%%%%%%%%%%%%%%%%%%%%%%%%%     Case 2 a    %%%%%%%%%%%%%%%%%%%%%%%%%%%%%%%%%%%%%%%%%%%%%%%%%%%%%%%%%%%%%%%%%%%%%%%%%%%%%%%%%%%%%%
%%%%%%%%%%%%%%%%%%%%%%%%%%%%%%%%%%%%%%%%%%%%%%%%%%%%%%%%%%%%%%%%%%%%%%%%%%%%%%%%%%%%%%%%%%%%%%%%%%%%%%%%%%%%%%%%%%%%%%%%%%%%%%%%%%%%%%%%%%%%%%%%%%%%%%%%

\textbf{Case 2a}: Consider $kC_{10}$  to be an arbitrary snake with a string as indicated in the condition of Case 2 with $k \geq 4$ and $k$ even. The proof is similar to the proof of Case 1a.% So, by the induction hypothesis, there is a near graceful labelling with vertex labels that are a subset of $\left[0,10k-9\right] \setminus \left\{2,10k-10\right\},$ and the edge labels are exactly $\left[1,10k-9\right] \setminus \left\{10k-10\right\}$, with $0$ in the last cycle in any position.

We label $kC_{10}$ as follows. For the first $k-1$ cycles, use the labelling of $\left(k-1\right)C_{10}$ obtained by induction and then add $4$ to each vertex label. Thus, the vertices have been labelled from the set $\left[4,10k-5\right] \setminus \left\{6,10k-6\right\}$ and by Lemma \ref{thirdlemma}, the edge labels are $\left[1,10k-9\right] \setminus \left\{10k-10\right\}$.

Apply the labelling $C_{10}^{c}$ to the final cycle, with the cut vertex receiving label $4$. Then there is a graceful labelling of $kC_{10}$ and $0$ in the $d=4$ position, from the cut vertex. (Note that a possible conflict occurs as $\left(10k-6\right)+4=10k-2$, however, in the labelling of the $(k-1)C_{10}$ no vertex is labelled $10k-2$, therefore we can use $C_{10}^{c}$ without any restriction).
%(Note that when we add $4$ in the labelling of $C_{10}^{a}$ we did not use $10k-6$, therefore we can use $C_{10}^{c}$ without any restriction).

We obtain a labelling with $0$ in the $d=3$ position by using the previously discussed labelling, ending with the $C_{10}^{c}$-labelling in the last cycle, then applying Lemma~\ref{third2lemma}.

\textbf{Case 2b}: Consider $kC_{10}$ to be an arbitrary snake with a string as indicated in the condition of Case 2 with $k \geq 3$ and $k$ odd. The proof is similar to the proof of Case 1b, using the labelling of $\left(k-1\right)C_{10}$ and adding $5$ to each vertex. Finally, apply the labelling $C_{10}^{a}$ to the final cycle, with the cut vertex receiving label $5$. Then we obtain a near graceful labelling of $kC_{10}$ and $0$ in the $d=4$ position, relative to cut vertex.

\textbf{Case 3}: $d_{i} \in \left\{2,3\right\}$ if $i$ is odd and $d_{i}=3$ if $i$ is even.

We prove that the $kC_{10}$ ($k \geq 1$) with string ($d_{1},d_{2},\ldots,d_{k-2}$) satisfying the condition of part three with $0$ in the $d=2$ or $d=3$ position in the last cycle and near graceful if $k$ is odd with $10k+1$ in the $d=3$ position in the last cycle.

For $k=1,$ use the labelling of $C_{10}^{d}$ in Table~\ref{c10c10} with $k=1$ which will make it near graceful. For $k=2,$ use $\left(12,6,13,8,\overline{16},5,14,10,11,9\right)$ with $(0$, $20$, $1$, $\overline{16}$, $3$, $17$, $7$, $19$, $2$, $18)$ to get $d=3$, and take this labelling with Lemma~\ref{third2lemma} to get $d=2$. The rest of the proof is similar to the proof of Case 1, but uses $C_{10}^{d}$ instead of $C_{10}^{a}$, and $C_{10}^{e}$ instead of $C_{10}^{b}$.

%%%%%%%%%%%%%%%%%%%%%%%%%%%%%%%%%%% C10 -- d=2,1 %%%%%%%%%%%%%%%%%%%%%%%%%%%%%%%%%%%%%%%%%%%%%%%%%%%%%%%%%%%%%%%%%%%%%%%%%%%%%%%%%%%%%%%%%%%%%%%%%%%%%%%%%%%%%%%%%%
\textbf{Case 4}: $d_{i} \in \left\{1,2\right\}$ if $i$ is odd and $d_{i}=2$ if $i$ is even.

We prove that $kC_{10}$ ($k \geq 1$) with string ($d_{1},d_{2},\ldots,d_{k-2}$) satisfying the condition of part four with $0$ in the $d=1$ or $d=2$ position in the last cycle and near graceful if $k$ is odd with $10k+1$ in the $d=1$ position in the last cycle.

For $k=1,$ use the labelling of $C_{10}^{f}$ in Table~\ref{c10c10} with $k=1$ which will make it near graceful. For $k=2,$ use $\left(8,14,7,12,\overline{4},15,6,10,9,11\right)$ with $(0$, $20$, $2$, $18$, $1$, $13$, $3$, $17$, $\overline{4}$, $19)$ to get $d=2$, and take this labelling with Lemma~\ref{third2lemma} to get $d=1$. The rest of the proof is similar to the proof of Case 2, but uses $C_{10}^{f}$ instead of $C_{10}^{a}$, and $C_{10}^{g}$ instead of $C_{10}^{c}$.\end{pf}
In Theorem~\ref{c14new} we prove (near) graceful labellings exist for a $kC_{14}$ for particular strings.
In Table~\ref{c1414} we see two labellings of $C_{14}.$ The labelling $C_{14}^{a}$ uses edge labels $\left[14k-13,14k+1\right] \setminus \left\{14k\right\}$. The labelling $C_{14}^{b}$ uses edge labels $\left[14k-14,14k\right] \setminus \left\{14k-13\right\}$.
%In Table~\ref{c1414} we see two labellings of $C_{14}.$ The labellings $G$ use vertex labels $\left\{0,1,2,4,5,6,7,14k-1,14k-2,14k-3,14k-4,14k-5,14k-6,14k+1\right\}$ and edge labels $\left\{14k+1,14k-1,14k-2,14k-3,14k-4,14k-5,14k-6,$\\$14k-7,14k-8,14k-9,14k-10,14k-11,14k-12,14k-13\right\}$}.\\ The labellings $H$ use vertex labels $\left\{0,1,2,3,4,5,6,14k,14k-1,14k-2,$\\$14k-3,14k-4,14k-5,14k-10\right\}$}, and edge labels $\left\{14k,14k-1,14k-2,$\\$14k-3,14k-4,14k-5,14k-6,14k-7,14k-8,14k-9,14k-10,14k-11,14k-12,14k-14}$.
%%%%%%%%%%%%%%%%%%%%%%%%%%%%%%%%%%%%%%%%%%%%%%%%%%%%%%%%%%%%%%%%%%%%%%%%%%%%%%%%%%%%%%%%%%%%%%%%%%%%%%%%%%%%%%%%%%%%%%%%%%%%%%%%%%%%%%%%%%%%%%%%%%%%%%%%
\begin{table}[H]
\large
\begin{center}\setlength\extrarowheight{7pt}\renewcommand{\arraystretch}{1.0}
\scalebox{0.70}{
\begin{tabular}{|c|c|c|}
\hline
 & Labelling    \\
\hline
$C_{14}^{a}$ &\ $\left(0,14k+1,4,14k-2,5,14k-4,6,14k-5,\overline{7},14k-6,2,14k-3,1,14k-1\right)$    \\  % it was G
\hline
$C_{14}^{b}$ &\ $\left(0,14k,4,14k-10,2,14k-3,3,14k-5,\overline{6},14k-4,5,14k-2,1,14k-1\right)$        \\ % it Was H
\hline

\end{tabular}}
\caption{Some useful labellings of $C_{14}$.} \label{c1414}
\end{center}
\end{table}
\begin{theorem}\label{c14new}
If $k \geq 1$ and $d_{i} \in \left\{6,7\right\}$ if $i$ is odd and $d_{i}=7$ if $i$ is even, for $1 \leq i \leq k-2$, then $kC_{14}$ with string ($d_{1},d_{2},\ldots,d_{k-2}$) is graceful if $k$ is even and near graceful if $k$ is odd.\end{theorem}
\begin{pf}The proof is similar to the proof of Theorem~\ref{ac6}.

As in the proof of Theorem~\ref{ac6}, we prove that if $k \geq 1$ and $d_{i} \in \left\{6,7\right\}$ if $i$ is odd and $d_{i}=7$ if $i$ is even, for $1 \leq i \leq k-2$, then the $kC_{14}$ with string ($d_{1},d_{2},\ldots,d_{k-2}$) is graceful if $k$ is even with $0$ in the $d=6$ or $d=7$ position in the last cycle and near graceful if $k$ is odd, with $14k+1$ in the $d=7$ position in the last cycle.% If $0$ is in the last cycle of a $kC_{14},$ so up to symmetry its position is uniquely determined by the distance from the last cut-vertex.% These distances, $d \in \left[1,7\right]$.

We proceed by induction on $k.$ For $k=1,$ use the labelling of $C_{14}^{a}$ in Table~\ref{c1414} which will make it near graceful. For $k=2$ use $(21$, $8$, $20$, $10$, $19$, $13$, $14$, $12$, $15$, $11$, $16$, $9$, $17$, $\overline{6})$ and $(0$, $28$, $4$, $18$, $2$, $25$, $3$, $23$, $\overline{6}$, $24$, $5$, $26$, $1$, $27)$ to get $d=6$, and take this labelling with Lemma~\ref{third2lemma} to get $d=7$.% by Theorem~\ref{oddcycle}.\\

\textbf{Case 1}: Consider $kC_{14}$ to be an arbitrary snake with $0$ in the $d=6$ or $d=7$ position in the last cycle with $k \geq 4$ and $k$ even. % Let $G$ the $(k-1)C_{14}$ obtained by deleting a last cycle from this $kC_{14}$. By the induction hypothesis, there is a near graceful labelling of $G$ with a $0$ on the vertex distance $d_{k-2}$ from the previous cut vertex.
The labelling of $(k-1)C_{14}$ by the induction hypothesis has vertex labels that are a subset of $\left[0,14k-15\right] \setminus \left\{3,14k-14\right\}$ and the edge labels are exactly $\left[1,14k+1\right] \setminus \left\{14k\right\}$ with $14k-13$ in the vertex distance $d_{k-2}$ from the previous cut vertex.

We label $kC_{14}$ as follows. Subtract each vertex label from $14k-7$ for the first $k-1$ cycles of the $(k-1)C_{14}$, so that the last cut vertex receives label $6$. Then, apply the labelling $C_{14}^{b}$ to the final cycle, with the cut vertex receiving label $6$. Then this labelling of $kC_{14}$ has all vertices labelled from $\left[0,14k\right]$ and the edge labels are exactly $\left[1,14k\right]$, with no repeated vertex or edge labels. Hence, we obtain a graceful labelling of $kC_{14}$ and $0$ in the $d=6$ position.

By using the previously discussed labelling, ending with the $C_{14}^{b}$-labelling in the last cycle, then applying Lemma~\ref{third2lemma} we obtain a labelling with $0$ in the $d=7$ position.
%%%%%%%%%%%%%%%%%%%%%%%%%%%%%%%%%%%%%%%%%%%%%%%%%%%%%%%%%%%%%%%%%%%%%%%%%%%%%%%%%%%%%%%%%%%%%%%%%%%%%%%%%%%%%%%%%%%%%%%%%%%%%%%%%%%%%%%%%%%%%%%%%%%%%%%%%%%%%%%%%%%%%%%%%%%%%%%%%%%%%%%%%%%%%%%%%%%%%%%%%%%     Case 2    %%%%%%%%%%%%%%%%%%%%%%%%%%%%%%%%%%%%%%%%%%%%%%%%%%%%%%%%%%%%%%%%%%%%%%%%%%%%%%%%%%%%%%
%%%%%%%%%%%%%%%%%%%%%%%%%%%%%%%%%%%%%%%%%%%%%%%%%%%%%%%%%%%%%%%%%%%%%%%%%%%%%%%%%%%%%%%%%%%%%%%%%%%%%%%%%%%%%%%%%%%%%%%%%%%%%%%%%%%%%%%%%%%%%%%%%%%%%%%%

\textbf{Case 2}: Consider $kC_{14}$ to be an arbitrary snake with $0$ in the $d=6$ position in the last cycle with $k \geq 3$ and $k$ odd. The proof is similar to the proof of Case 1. Add $7$ to each vertex label instead of subtracting each vertex label from $14k-7$ and apply the labelling $C_{14}^{a}$ to the final cycle instead of $C_{14}^{b}$, with the cut vertex receiving label $7$ instead of $6$. Then this labelling of $kC_{14}$ has all vertices labelled from $\left[0,14k+1\right] \setminus \left\{3,14k\right\}$, and the edge labels are exactly $\left[1,14k+1\right] \setminus \left\{14k\right\}$, with no repeated vertex or edge label. Thus, we obtain a near graceful labelling of $kC_{14}$ and $14k+1$ in the $d=7$ position.\end{pf}
%%%%%%%%%%%%%%%%%%%%%%%%%%%%%%%%%%%%%%%%%%%%%%%%%%%%%%%%%%%%%%%%%%%%%%%%%%%%%%%%%%%%%%%%%%%%%%%%%%%%%%%%%%%%%%%%%%%%%%%%%%%%%%%%%%%%%%%%%%%%%%%%%%%%%%%%%%%%%%%%
Recall that a $kC_{t}$ is linear if all entries in its string are equal to $\left\lfloor \frac{t}{2}\right\rfloor$. So, based on the results of Theorems~\ref{c10new} and~\ref{c14new} we now state a corollary for linear $kC_{10}$ and $kC_{14}$, following the style of Theorem~\ref{t2}.   
\begin{corollary}\label{cc14}If $k \geq 1$ then the linear $kC_{10}$ and linear $kC_{14}$ are graceful if $k$ is even and nearly graceful if $k$ is odd.\end{corollary}
\section{Discussion}
In this paper we (near) gracefully labelled several type of snakes. In section~\ref{ss2}, we presented a new sufficient condition which when satisfied shows there is a graceful labelling of a $kC_{4n}$ for any string. By using a complete $s,t$-useful cycle set we proved that if there is a complete $s,4n$-useful cycle set with $s \geq 1$, then there exists a graceful labelling of any $kC_{4n}$. We used the results in \cite{rosa} with our results for $kC_{4n}$ and proved that a graceful labelling exists for particular $kC_{4n}$ with string $\left(d_{1},d_{2},\ldots,d_{k-2}\right)$, where $d_{i} \in \left\{2,3,4\right\}$. Expanding these results for any $n$ and $d$ is possible but hard to apply for large $n$. We extended our main result to the case of cyclic snakes with cycles of varying sizes. Further, we extended the results in Theorems~\ref{t2} to~\ref{t5} on (near) gracefully labelled $kC_n$ where $n = 6,8,12,16,20,24$ for all possible snakes.
%In Section~\ref{ss3}, as we discussed we have to find two complete useful cycle sets for $4n+2$ because we are trying to change a graceful labelling to a near-graceful one, or the reverse. If we have these complete cycles then we can prove the (near) graceful labelling for any $kC_{4n+2}$. For $kC_{6}$ we found the complete cycle set and we gave the spectrum solution for $kC_{6}$ in Theorem~\ref{ac6}. For $kC_{10}$ and $kC_{14}$ we presented new results on the (near) graceful labelling of $kC_n$ where $n = 10,14$ and $k > 1$. But these results not complete because we are unable to find the complete cycles sets for them. The result of Section~\ref{ss3} leads us to ask the following question:
%\textbf{Question:} Can we gracefully label any $kC_{4n}$-snake with $n \geq 1$? we conjecture that the answer is yes.

As we discussed in Section~\ref{ss3}, new approaches must be found to gracefully label $kC_{4n+2}$ snakes, even for fixed $n$. Our collections of ad hoc methods work to give classes for fixed $n$, but do not seem to generalize, even for ``nice'' subfamilies, such as linear snakes. Thus we pose the following open question.

\textbf{Question:} Can we (near) gracefully label every $kC_{m}$ with $k \geq 1$ and $m\equiv  2\ ($mod$\ 4)$?

In fact the technique we used is more general than indicated in our theorems. Suppose we have a gracefully labelled bipartite graph $G=K_{3,4}$ as in Figure~\ref{go}. If we add $3$ to each vertex label and use the cycle $H=C_{8}^{2}$ from Table~\ref{c4n4n} then we obtain a new gracefully labelled graph as in Figure~\ref{gn}. Thus we can in several cases gracefully label new graphs.
 \begin{theorem}If $G$ is graceful and $H$ is a $kC_{4n}$ with $1 \leq n \leq 6$, then the graph $GH^{\ast}$ obtained by identifying any vertex in $G$ that can be labelled $0$ in some graceful labelling with any vertex in the first cycle of $H$ is graceful.\end{theorem}
%then we obtain a new gracefully labelled graph $GH^{\ast}$, where $GH^{\ast}$ is a graph obtained by by identify the $2n-1$ vertices in the $G$ and $H$ labelling. .
%%%%%%%%%%%%%%%%%%%%%%%%%%%%%%%%%%%%%%%%%%%%%%%%%%%%%%%%%%%%%%%%%%%%%%%%%%%%%%%%%%%%%%%%%%%%%%%%%%%%%%%%%%%%%%%%%%%%%%%%%%%%%%%%%%%%%%%%%%%%%%%%%%%%%%%%%%%%%%%%%%%%%%%%
\begin{figure}[H]
\begin{center}\setlength\extrarowheight{7pt}
  \begin{minipage}[b]{0.4\textwidth}
\begin{tikzpicture}[c/.style={circle,inner sep=1pt,fill}]
 \path foreach \x [count=\y] in {12,8,4,3,2,1,0}
  {({0.8*sign(\y-3.5)},{4-Mod(\y,4)-0.25*sign(\y-3.5)})
    node[c,label={90-90*sign(\y-3.5)}:{$\x$}](\x){}};
  \path foreach \x in {12,8,4}
   {foreach \y in {3,2,1,0} {(\x) edge (\y)}};
\end{tikzpicture}
		\caption{Gracefully labelled $K_{3,4}$.}\label{go}
		  \end{minipage}
\quad
 %\hfill 	
		%\end{center}
	%	\end{figure}
%%%%%%%%%%%%%%%%%%%%%%%%%%%%%%%%%%%%%%%%%%%%%%%%%%%%%%%%%%%%%%%%%%%%%%%%%%%%%
%\begin{figure}[H]
%\begin{center}\setlength\extrarowheight{7pt}
  \begin{minipage}[b]{0.39\textwidth}
\begin{tikzpicture}[c/.style={circle,inner sep=1pt,fill}]
 \path foreach \x [count=\y] in {15,11,7,6,5,4,3}
  {({0.8*sign(\y-3.5)},{4-Mod(\y,4)-0.25*sign(\y-3.5)})
    node[c,label={90-90*sign(\y-3.5)}:{$\x$}](\x){}};
  \path foreach \x in {15,11,7}
   {foreach \y in {3,4,5,6} {(\x) edge (\y)}};
	
	\draw [thick, -] (1.75,-0.25) -- (2.75,-0.25) -- (3.75,-0.25) -- (4.75,0.75) -- (3.75,1.75) -- (2.75,1.75) -- (1.75,1.75) -- (0.75,0.75)  -- (1.75,-0.25);
\draw[fill] (1.75,-0.25) circle  [radius=.9pt];\draw[fill] (2.75,-0.25) circle  [radius=.9pt];\draw[fill] (3.75,-0.25) circle  [radius=.9pt];\draw[fill] (4.75,0.75) circle  [radius=.9pt];\draw[fill] (3.75,1.75) circle  [radius=.9pt];\draw[fill] (2.75,1.75) circle  [radius=.9pt];\draw[fill] (1.75,1.75) circle  [radius=.99pt];\draw[fill] (0.75,0.75) circle  [radius=.9pt];%;\draw[fill] (0.75,0.75) circle  [radius=.9pt];

\node [below] at (1.75,-0.25) { $20$};
\node [below] at (2.75,-0.25) { $0$};
\node [below] at (3.75,-0.25) { $19$};
\node [right] at (4.75,0.75) { $1$};
\node [above] at (3.75,1.75) { $17$};
\node [above] at (2.75,1.75) { $2$};
\node [above] at (1.75,1.75) { $16$};
\node [left] at (0.75,0.75) { $$};
\end{tikzpicture}
		\caption{Gracefully labelled $GH^{\ast}$.}\label{gn}
		\end{minipage}
		\end{center}
		\end{figure}

\end{document}